\documentclass[12pt,draftclsnofoot,onecolumn]{IEEEtran}

\usepackage[utf8]{inputenc}
\usepackage{xcolor}

\usepackage{mathptmx} 
\usepackage{amsmath,amssymb} 
\usepackage{dsfont}
\usepackage{verbatim}
\usepackage{multicol,longtable}
\usepackage{epsfig}
\usepackage{calc}
\usepackage{subcaption}
\usepackage{enumerate}
\usepackage{bm}
\usepackage{color}
\usepackage{url}
\usepackage{algorithm}
\usepackage[noend]{algpseudocode}

\newtheorem{remark}{Remark}[section]

\begin{document}

\title{2D Density Control of Micro-Particles using Kernel Density Estimation}


\author{ Ion Matei, Johan de Kleer and Maksym Zhenirovskyy
\thanks{Ion Matei, Johan de Kleer and Maksym Zhenirovskyy are with the Intelligent Systems Laboratory, at Palo Alto Research Center (emails: {\tt imatei@parc.com, dekleer@parc.com, mazhenir@parc.com}).}
}

\maketitle
\begin{abstract}
We address the problem of 2D particle density control. The particles are immersed in dielectric fluid and acted upon by manipulating an electric field. The electric field is controlled by an array of electrodes and used to bring the particle density to a desired pattern using dielectrophoretic forces. We use a lumped, 2D, capacitive-based, nonlinear model describing the motion of a particle. The spatial dependency of the capacitances is estimated using electrostatic COMSOL simulations. We formulate an optimal control problem, where the loss function is defined in terms of the error between the particle density at some final time and a target density. We use a kernel density estimator (KDE) as a proxy for the true particle density. The KDE is computed using the particle positions that are changed by varying the electrode potentials. We showcase our approach through numerical simulations, where we demonstrate how the particle positions and the electrode potentials vary when shaping the particle positions from a uniform to a Gaussian distribution.
\end{abstract}

\section{Introduction}
\begin{figure}[!htp]
  \begin{center}
    \includegraphics[width=0.9\textwidth]{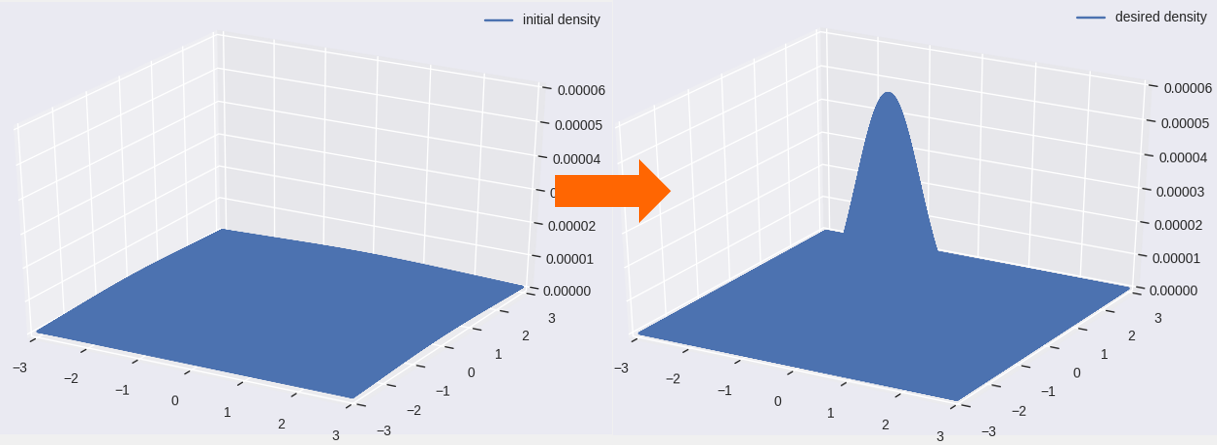}
    \end{center}
  \caption{Control objective: given an initial distribution of particles apply a sequence of electrode potentials to shape the particle distribution in a desired pattern.}
  \label{fig:03261834}
\end{figure}

We aim to design and build a printer system for assembling micro-particles into engineered patterns. Micro-particles are submerged into dielectric fluid and their positions are controlled by manipulating the electric potential of a 2D array of electrodes. The particle positions are tracked by a high-speed camera and the control signals are generated by projecting images on photo-sensitive transistors attached to the electrodes. In
\cite{7994140,8709705,7963173,9303808} we introduced control algorithms that act on individual particles. However, controlling and tracking of a large number of particles, at an individual level, is computationally expensive. It is more advantageous to control  simultaneously large number of particles and shape them into a desired density. This is the problem we are addressing in this paper. At a conceptual level, the control objective is graphically depicted in Figure \ref{fig:03261834}: we start with an initial particle density (e.g., uniform distribution) and we would like to converge to a target particle density (e.g., multivariate Gaussian distribution), over some time horizon, by varying the electrode potentials.
There is a close connection between our problem and the optimal control of the Liouville  equation. In fact, as shown later in paper, our formulation can be mapped to an optimal control problem where the dynamical constraint is the Liouville partial differential equation (PDE) in terms of a probability density function. Controllability of the Liouville equation together with optimal control of its moments for some special cases (e.g., linear case) are discussed in \cite{Brockett2012}. An analysis of problems of optimal control of ensembles governed by
the Liouville equation is done in \cite{Bartsch2019ATI}, where the results apply to particular classes of problems (e.g., Liouville equation with unbounded drift
function with linear and bilinear control mechanisms), and classes of cost functionals. In \cite{8618429}, the authors introduce a dynamic output feedback control of the Liouville equations, applied to SISO discrete-time linear systems. Optimal control of the Liouville equation is part of a larger class control problems having PDEs as dynamical constraints.  The type of PDE typically determines the methods for solving their related analysis and synthesis problems. For parabolic PDE systems, their dominant dynamics can be characterized by finite-dimensional, ordinary differential equation (ODE) systems, generated through Galerkin projection, for example. The control approach in this case is based on using the ODEs for controller design \cite{10.1115/1.1451164,1239777}. In the case of hyperbolic PDE systems, an infinite number of modes are needed to represent their dynamics. As a consequence, control schemes are designed by taking into account the spatial dimension, as well \cite{10.1016/j.automatica.2009.02.017,doi:10.1021/ie0341404,SHANG2000533}. Another control application under PDE dynamical constraints is traffic control. For example, linearized versions of the Aw-Rascle-Zhang PDE model can be  used to determine an output feedback control law \cite{yu:hal-02732932} that can be applied to traffic control.

While solving an optimal control problem in terms of the Liouville equation is feasible, such an approach is numerically complex. The complexity stems from the need to generate a discrete representation for the Liouville PDE, either through a finite elements approach or using spectral methods. In this paper we use a particle-based approach to approximate the particle density. In particular, we use a kernel density estimator (KDE) as a proxy for the particle density and solve the optimal control problem in terms of this quantity. Estimating the initial density would be needed even when using the Liouville equation, since this estimate is used as initial value.  Our objective is to compute a sequence of electrode electric potentials so that an initial particle distribution is shaped into a target distribution after applying this sequence over time.  We define the optimal control cost function in terms of the $L_2$ metric used to compute the error between the particle density at the end of a time horizon and a target density. Some authors \cite{DevroyeGL96} make a strong case for using an $L_1$ metric to compute the error since it is transformation invariant. However, it is more difficult to deal with from a numerical optimization perspective. Another possible loss function is the Kullback-Leibler loss, but it is not the best choice for non-parametric densities since it is completely dominated by the tails of the densities. The KDE depends on the predicted trajectories of a set of particles, where the trajectory of a single particle is determined by a lumped, 2D, capacitive-based, nonlinear model describing its motion. We assume that there is no interaction between particles and that their motion is completely determined by the electric field generated by the array of electrodes. The electric field induces an accumulation of potential energy at the particles. We use automatic differentiation (AD) enabled by {\tt jax} \cite{jax2018github} to compute the forces (i.e., the gradient of the potential energy) that act on the particles, and the gradients of the loss and constraint functions that are passed to the optimization algorithm.

\emph{{Notations}:} We denote scalars, vectors and random variables by Italic symbols, bold Italic symbols, and capital Italic symbols, respectively.  Let $f(\boldsymbol{x})$ be a multi-variable map, where $\boldsymbol{x} = (x_i)$ is a vector of scalars.  We denote by $\nabla f(\boldsymbol{x})$ the gradient of $f$, and by $\frac{\partial f}{\partial x_i}(\boldsymbol{x})$ the partial derivative of $f$, with respect to $x_i$. 
For a vector valued function $\boldsymbol{F}(\boldsymbol{x})$, $\frac{\partial \boldsymbol{F}}{\partial \boldsymbol{x}}$ denotes the Jacobian of $\boldsymbol{F}$. For a function $f(\boldsymbol{x})$, with $\boldsymbol{x} = (x_1, x_2)$, $f(x_1, x_2)$ is an equivalent notation. The divergence operator applied to a vector-valued function $\boldsymbol{F}(\boldsymbol{x})$ is denoted by $\nabla \cdot \boldsymbol{F}$. For a matrix $A$, $|A|$ denotes its determinant.

\emph{{Paper structure}:} In Section \ref{sec:Particle motion control setup} we introduce the experimental setup for controlling the micro-particles, and the single particle model of motion. 
We formulate the particle density control problem in Section \ref{sec:Feedback control} and showcase our approach through simulation results in Section \ref{sec:Results}.

\section{Particle motion control setup}
\label{sec:Particle motion control setup}
We first introduce the experimental setup for controlling the particle density. Second, we describe the dynamical model for a single particle when actuated by an array of electrodes. Third, we generate a dynamical model for a single particle under the actuation of an arbitrarily large number electrodes.

\subsection{Experimental setup}
\begin{figure}[!htp]
  \begin{center}
    \includegraphics[width=0.9\textwidth]{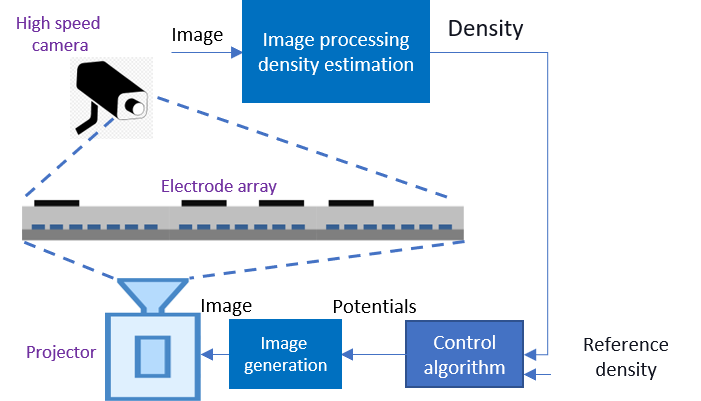}
    \end{center}
  \caption{Diagram of the experimental setup: a high-speed camera tracks the particle density that is transmitted to the control module. A control module generates electrode potentials that minimizes the error between a target density and the particle density and the end of a time horizon. The control inputs (electrodes potential) are converted into images that are projected on the array. The light produced by the images powers the photo-transistors attached to the electrodes.}
  \label{fig:03271129}
\end{figure}

The experimental setup is shown in Figure \ref{fig:03271129}. The system has three hardware devices and three software modules. The hardware devices include: a high-speed camera for tracking the particle densities, an array of electrodes to generate a dynamic potential energy landscape for manipulating objects with dielectrophoretic (DEP) forces, and a video projector to actuate the array based on projected images. The software modules include: a module for image processing that estimates the particle density, a control module that compares the target particle density with the density at the end of the control horizon, and generates input signals to minimize the error between them, and an image generation module that maps the control inputs to images that are projected on the array. The projected images activate or deactivate electrodes, as indicated by the control inputs. 

\subsection{Single particle dynamical model}
\label{lab:particle_model}
\begin{figure}[ht!]
\begin{center}
\includegraphics[width=0.8\textwidth]{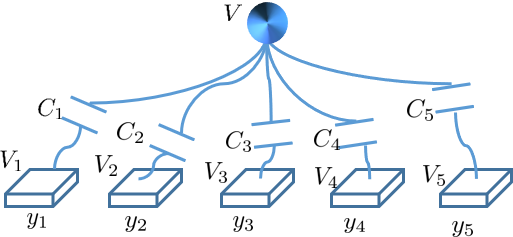}
\end{center}
\caption{Capacitive-based model describing the interaction between particle and electrodes situated at positions $y_i$.}
\label{fig:11181642}
\end{figure}

\begin{figure*}[ht!]
\begin{center}
\includegraphics[width=0.8\textwidth]{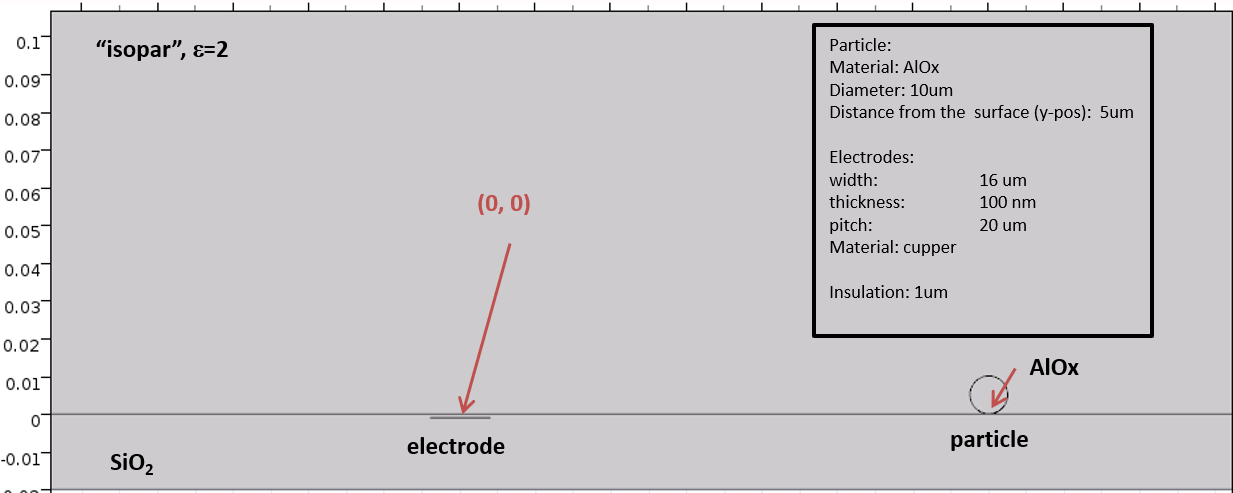}
\end{center}
\caption{Electrostatic COMSOL model of two conductors (spherical particle and electrode) for capacitance computation. The electrode's position is $(0,0)$. }
\label{fig:11181648}
\end{figure*}

\begin{figure}[!htp]
  \begin{center}
    \includegraphics[width=0.7\textwidth]{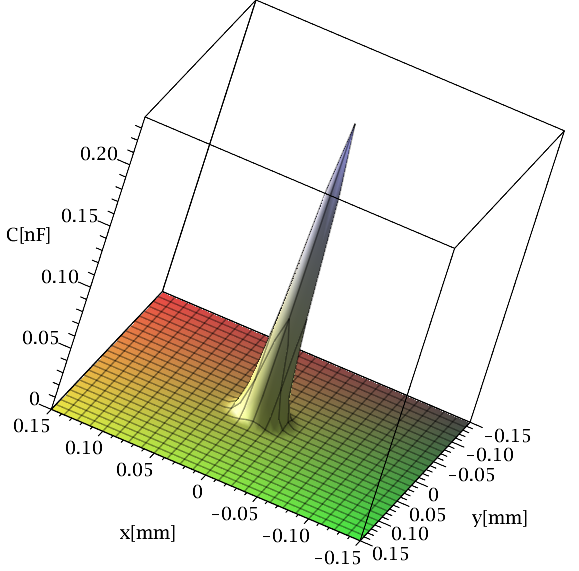}
    \end{center}
  \caption{2D capacitance function determined from the 1D model using the symmetry property.}
  \label{fig:2d_capacitance}
\end{figure}

In what follows we described a 2D model for the particle motion under the effect of the potential field induced by the electrode array. The model is for one particle only and neglects possible interactions when particles get close to each other. By applying electric potentials to the electrodes, we generate DEP forces that act on the particles. A viscous drag force proportional to the velocity\footnote{The drag force is proportional to the velocity in non-turbulent flows, that is, when the Reynolds number is small.} opposes the particle's motion. Due to the negligible mass of the particle, the acceleration can be neglected. It follows that the particle dynamical model can be described by:
\begin{equation}
\label{eq:11191625}
\mu \dot{\boldsymbol{x}} = F(\boldsymbol{x}, \boldsymbol{V}),
\end{equation}
where $\boldsymbol{x}$ denotes the 2D particle position measured at its center of mass, $\boldsymbol{V} = (V_{i,j})$ are the electrode electric potentials,  $\mu$ is the fluid dependent viscous coefficient, and $F(\boldsymbol{x},\boldsymbol{V})$ is the vector of forces acting on the particle. The indices $(i,j)$ are associated to 2D electrode positions $\boldsymbol{y}_{i,j}$. We express the forces $F(\boldsymbol{x},\boldsymbol{V})$ as a function of the potential energy of the particle. We compute the potential energy using a capacitive-based electrical circuit that lumps the interaction between the electrodes and the particle. Such a circuit is shown in Figure \ref{fig:11181642}, where only one row with five electrodes of the array is depicted. The particle and the electrodes act as metal plates; hence the capacitances of these capacitors are dependent on the particle position. As expected, the maximum values are attained when the particle's position maximizes the overlap with the electrodes. To simplify the analysis, we limited our analysis to low frequency region only, where the dielectric constant is not frequency dependent.
The vector of forces $F$ can be formally expressed as $F(\boldsymbol{x},\boldsymbol{V})=\nabla U(\boldsymbol{x},\boldsymbol{V})$ , where $U(\boldsymbol{x},\boldsymbol{V})$ is the potential energy of the particle, given by:
\begin{equation}
\label{equ:11191636}
U(\boldsymbol{x},\boldsymbol{V}) = \frac{1}{2}\sum_{i,j=1}^NC_{i,j}(\boldsymbol{x})\left[V_{i,j}-v(\boldsymbol{x},\boldsymbol{V})\right]^2,
\end{equation}
where $C_{i,j}(\boldsymbol{x})$ is the capacitance between the particle at position $\boldsymbol{x}$ and electrode $(i,j)$, $V_{i,j}$ is the electric potential of electrode $(i,j)$, $v(\boldsymbol{x})$ is the electric potential of the particle, and $N$ is the number of actuated electrodes. If there are electrophoretic effects on the particle, we can readily extend the potential energy to include such effects. We compute the particle potential in terms of the electrode potentials, by solving for the voltages and currents in the electrical circuit shown in Figure \ref{fig:11181642}. In particular, the steady state particle potential is given by
\begin{equation}
\label{equ:11191701}
v(\boldsymbol{x},\boldsymbol{V}) = \frac{1}{\sum_{i,j=1}^NC_{i,j}(\boldsymbol{x})}\sum_{i,j=1}^NC_{i,j}(\boldsymbol{x})V_{i,j}.
\end{equation}
Feedback control design requires explicit expressions for the capacitances between the particle and electrodes. We learn the capacitance model using COMSOL simulations. In what follows, we describe through an example the process for learning the capacitances. This process can be repeated for other types of particles, geometries and material properties.  For symmetric particles (e.g., beads), we estimate the capacitances by simulating a 2D electrostatic COMSOL model. This implies that the capacitance function is of the form $C_{i,j}(\boldsymbol{x}) = C(\|\boldsymbol{x}-\boldsymbol{y}_{i,j}\|)$, where $\boldsymbol{x}=(x_1,x_2)$ denotes the particle 2D position, and $\boldsymbol{y}_{i,j}$ is the fixed position of electrode $i,j$.  As shown in Figure \ref{fig:11181648}, in the COMSOL model, a 16 $\mu m$ width and 100 $nm$ thickness copper plate,  and a 10 $\mu m$ diameter,  aluminium oxide (AlOx) spherical object are surrounded by a dielectric  with the same properties as the isopar-M solution.
The quasi-static models are computed in form
of electromagnetic simulations using partial
differential equations, where we use ground boundary (zero potential) as the boundary condition. The capacitance matrix entries
are computed from the charges that result on each conductor when an electric potential is applied to one of them and the other is set to ground. 
The COMSOL electrostatic model has as parameters, the diameter of the sphere, the electrode dimensions, the dielectric fluid constant ($\varepsilon=2$), the positions and material of the particle and the electrode. We fix the particle height at 5 $\mu m$ and generate simulation results by varying its position on the $x_1$-axis over the interval $[-1 mm, 1 mm]$. Note that due to the size of the particle versus the size of the electrodes, fringe effects (electric field distortions at the edges) are significant. The simulation results generate capacitances between the electrode and the particle for all considered positions. We parameterized the capacitance function using error functions: $C(\xi) = \sum_{i=1}^m a_i\left[ \Phi\left(\frac{\xi+\delta}{c_i}\right)-\Phi\left(\frac{\xi-\delta}{c_i}\right)\right]$, where $\Phi(\xi) = \frac{1}{\sqrt{\pi}}\int_{-\xi}^\xi e^{-t^2}dt$ is the error function, $\xi$ is the distance between the center of the particle and the electrode center assumed at the origin, $a_i$ and $c_i>0$ are scalars, and $\delta$ is half of the electrode pitch (10 $\mu m$ in our example).
\begin{figure}[!htp]
  \begin{center}
    \includegraphics[width=0.7\textwidth]{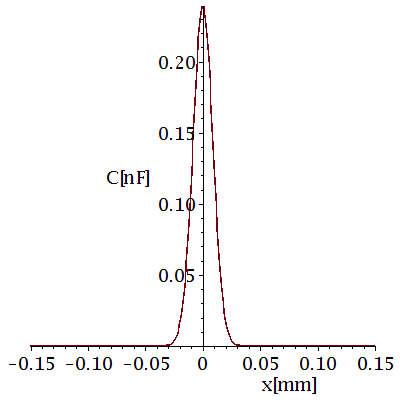}
    \end{center}
  \caption{Capacitance function fitted with training data generated by COMSOL simulations.}
  \label{fig:capacitance}
\end{figure}
Figure \ref{fig:capacitance} depicts $C(\xi)$, the capacitance between the particle and the electrode as a function of the particle horizontal position,  where the numerical values were fitted on the error function parameterization. For symmetric particles (e.g., sphere shaped), we can map the 1D model to a 2D model using the transformation $\xi\rightarrow \sqrt{x_1^2+x_2^2}$, which results in a capacitance function
$C(x_1,x_2) = \sum_{i=1}^m a_i\left[ \Phi\left(\frac{\sqrt{x_1^2+x_2^2}+\delta}{c_i}\right)-\Phi\left(\frac{\sqrt{x_1^2+x_2^2}-\delta}{c_i}\right)\right]$. For the sphere shape particle, it was enough to consider only one term in the parameterization of the capacitance function, and the resulting 2D  shape  is depicted in Figure \ref{fig:2d_capacitance}.

\subsection{From discrete to continuous electrode actuation}
In this section, we demonstrate how  to transform the particle motion model from a discretized actuation mechanism (discrete set of electrodes) to a continuous one. We start with the 1D case. We represent the  particle electric potential  as $v(x, \boldsymbol{V}) = \frac{1}{\sum_{i=1}^NC_i(x)}\sum_{i=1}^NC_i(x)V_i$, and  interpret $v(x,\boldsymbol{V})$ as the expected value of a random function $\mathcal{V}(Y)|X=x$ over a discrete distribution $p_i(x) = C_i(x)/\sum_{i=1}^NC_i(x)$. We represent the probability mass function $p_i(x)$ as a conditional probability $p_i(x) = Pr(Y=y_i|X=x)$, and hence the particle potential can be expressed as $v(x) = E[\mathcal{V}(Y)|X=x]$, where $V_i=\mathcal{V}(y_i)$ is a function that reflects the electric potential at each point $y_i$. The discrete probability distribution can be seen as a discretization of a continuous probability distribution, i.e., $p_i(x) = \int_{y_i-\delta}^{y_i+\delta} f(y|x)dy$, where $\delta$ is half of the electrode pitch. The parameterization of the capacitance function in terms of the error functions tells us that the conditional probability density function (p.d.f.) is a mixture of Gaussian functions. For a sphere shaped particle, the mixture has only one term, and hence the capacitance is expressed as $C_i(x) = a\left[ \Phi\left(\frac{x-y_i+\delta}{c}\right)-\Phi\left(\frac{x-y_i-\delta}{c}\right)\right]=2a\int_{y_i-\delta}^{y_i+\delta}f(y|x)dy$, where $f(y|x) = \frac{1}{\sqrt{2\pi \sigma^2}}e^{-\frac{(y-x)^2}{2\sigma^2}}$, with $\sigma = \frac{c}{\sqrt{2}}$. Therefore, the particle potential in the continuous representation can be expressed as $\bar{v}(x)=E[\mathcal{V}(Y)|X=x]$, where the expectation is computed with respect to the conditional Gaussian distribution, $\mathcal{N}(x,\sigma)$, and $\mathcal{V}(Y)$ is a function that assigns an electric potential to each point $y$. The potential energy can now be represented as $U(x) = a E[(\mathcal{V}(Y)-\bar{v}(X))^2|X=x]$, and it follows that the 1D particle dynamics is given by the following partial differential equation
\begin{equation}
    \label{equ:19201120}
    \mu \dot{ x}= \frac{\partial U(x)}{\partial x},
\end{equation}
where $\bar{v}(x)=E[V(Y)|X=x]$.\\%
\emph{Extension to 2D case}: We denote by $\boldsymbol{x} = (x_1,x_2)$ the particle position and by $\boldsymbol{y}_{i,j}=(y_{1}^{i,j},y_{2}^{i,j})$ the position of electrode $(i,j)$. The particle dynamics becomes
\begin{equation}
    \label{equ:09251120}
    \mu \dot{\bm{x}} =\nabla U(\boldsymbol{x}, \boldsymbol{V}),
\end{equation}
where $\mu$ is the viscous coefficient, and $U(\boldsymbol{x}, \boldsymbol{V})$ denotes the particle's potential energy. As is the 1D case, the potential energy is given by
$$U(\boldsymbol{x}, \boldsymbol{V}) = \frac{1}{2}\sum_{i,j=1}^NC_{i,j}(\bm{x})\left[V_{i,j} - v(\bm{x})\right]^2,$$
where $v(\boldsymbol{x}, \boldsymbol{V})$ denotes the particle electric potential, $C_{i,j}(\bm{x})=C(\|\bm{x}-\bm{y}_{i,j}\|)$ represents the capacitance between the particle at the position $\boldsymbol{x}$ and electrode ${i,j}$ at position $\boldsymbol{y}_{i,j}$, and $V_{i,j}$ represents the potential of electrode $(i,j)$. Similar to the 1D case, we assume that the capacitance $C_{i,j}$ can be represented as the un-normalized discretization of a multi-variable Gaussian p.d.f., that is:
\begin{equation}\label{eq:03250923}
  C_{i,j}(\bm{x})=4a \int_{y_1^{i,j}-\delta/2}^{y_1^{i,j}+\delta/2}\int_{y_2^{i,j}-\delta/2}^{y_2^{i,j}+\delta/2}f(\bm{y}|\bm{x})d\bm{y},
\end{equation}
where the conditional density function $f(\bm{y}|\bm{x})$  is the multivariate Gaussian distribution $\mathcal{N}(\bm{X},\sigma^2I)$. The potential energy is similar to the 1D case and it is given by $U(\boldsymbol{x}) = a E[(\mathcal{V}(\boldsymbol{Y})-\bar{v}(\boldsymbol{X}))^2|\boldsymbol{X}=\boldsymbol{x}]$,
where $\mathcal{V}(\boldsymbol{y})$ is a map such that $\mathcal{V}(\boldsymbol{y}_{i,j}) = V_{i,j}$.
\begin{remark}
Using a mean-field approximation argument (see for instance \cite{Carrillo2010} in the context of the Cucker-Smale model) we can derive the Liouville equation that describes the evolution of the density of a large number of particles. We obtain a control input dependent PDE of the form:
\begin{equation}
\label{eqn:03272309}
\nabla \cdot \left[f(\boldsymbol{x},t) F(\boldsymbol{x},\boldsymbol{V}(t))\right]+\frac{\partial f(\boldsymbol{x}, t)}{\partial t} = 0,
\end{equation}
where $f(\boldsymbol{x},t)$ denotes the particle density, and $\boldsymbol{V}(t)$ are the control inputs (i.e., the electrode potentials). In this paper, we \emph{do not} use a  control approach based on a dynamics governed by the Liouville, but rather an optimal control approach where the dynamics is governed by a set of particles whose trajectories are used to approximate the particle density.
\end{remark}

\section{Feedback control}
\label{sec:Feedback control}
We formulate an optimal control problem to shape an initial density function $f_0(\boldsymbol{x}) = f(\boldsymbol{x},0)$ into a target density $f_d(\boldsymbol{x})$.
The ideal optimal control problem we would like to solve is: given a finite time horizon $[0, T]$, find the electric potentials $\boldsymbol{V}(t)=(V_{i,j}(t))$ for $t\in [0, T]$ such that $f(\boldsymbol{x},T)$ equals the target density $f_d(\boldsymbol{x})$. The main constraint of the optimization problem are the dynamics of the particle motion. We use a proxy for the particle density $f(\boldsymbol{x},T)$ given by a KDE. Let $\{\boldsymbol{x}^{(i)}\}_{i=1}^n$ be a set of $n$ particles. Then the KDE is given by
$$\hat{f}_H(\boldsymbol{x},t) = \frac{1}{n}\sum_{i=1}^n K_H\left(\boldsymbol{x} - \boldsymbol{x}_i\right),$$
where $H\in \mathds{R}^2$ is the symmetric, positive definite bandwidth matrix, $K_H(\boldsymbol{x}) = |H|^{-1/2}K(H^{-1/2}\boldsymbol{x})$, with $K$ being the kernel function. Examples of commonly used kernels include: boxcar, Gaussian, Epanechnikov or tricube. In this paper, we use the standard multivariate kernel: $K_{H} = \frac{1}{\sqrt{(2\pi)^2 |H|}}e^{-\frac{1}{2}\boldsymbol{x}^TH\boldsymbol{x}}$, due to its smoothness.

Let $\mathcal{X}\in \mathds{R}^2$ be a compact set that bounds the particle positions, and let $V_{max}$ be the maximum magnitude of the electric potentials.  We formulate the following optimal control problem:
\begin{IEEEeqnarray}{rcl}
\label{eqn:03272314}
\min_{\boldsymbol{V}(t),\boldsymbol{x}^{(i)}(t), t \in [0, T]} & & \int_{\mathcal{X}} (\hat{f}(\boldsymbol{x}, T)-f_d(\boldsymbol{x}))^2d\boldsymbol{x}\\
\nonumber
\textmd{s.t.: } &  & \dot{\boldsymbol{x}}^{(i)} = F(\boldsymbol{x}^{(i)}, \boldsymbol{V}(t)), \boldsymbol{x}^{(i)}(0)=\boldsymbol{x}^{(i)}_{0}\\
\nonumber
& & \hat{f}(\boldsymbol{x},T) = \frac{1}{n}\sum_{i=1}^n K_H\left(\boldsymbol{x} - \boldsymbol{x}^{(i)}(T)\right),\\
\nonumber
& & |\boldsymbol{V}(t)|\leq {V}_{max}, \forall t\in [0, T],\\
\nonumber
& & \boldsymbol{x}^{(i)}(t)\in\mathcal{X}, \forall t\in [0, T],
\end{IEEEeqnarray}
where $\boldsymbol{x}^{(i)}_{0}$ are the initial particle positions.

To convert the problem (\ref{eqn:03272314}) into a format amenable to numerical optimization, we need to discretize both the time and space, and generate approximations of the loss function and of the particle dynamics. For the spatial discretization, we can employ an uniform mesh with  cells that are centered at $(\boldsymbol{x}_{i,j})$, each cell having area $\Delta x^2$. In addition, we discretize the time using a sampling period $\Delta t$, resulting in a sequence of samples $\{t_k\}$. With these discretization schemes, and employing a trapezoidal rule to approximate the particle dynamics, the optimization problem (\ref{eqn:03272314}) becomes:
\begin{IEEEeqnarray}{rcl}
\label{eqn:03242113}
\min_{\boldsymbol{V}(t_k),\boldsymbol{x}^{(i)}(t_k)} & & \sum_{i,j} (\hat{f}(\boldsymbol{x}_{i,j}, T)-f_d(\boldsymbol{x}_{i,j}))^2\\
\nonumber
\textmd{s.t.: } &  & \boldsymbol{x}_{t_{k+1}}^{(i)} = \boldsymbol{x}_{t_{k}}^{(i)} + \frac{\Delta t}{2}\left[F(\boldsymbol{x}^{(i)}(t_{k+1}), \boldsymbol{V}(t_{k+1}))\right.\\
\nonumber
& &\left. +F(\boldsymbol{x}^{(i)}(t_{k}), \boldsymbol{V}(t_{k}))\right], \boldsymbol{x}^{(i)}(0)=\boldsymbol{x}^{(i)}_{0}\\
\nonumber
& & \hat{f}(\boldsymbol{x}_{i,j},T) = \frac{1}{h}\sum_{i=1}^n K_H\left(\boldsymbol{x}_{i,j} - \boldsymbol{x}^{(i)}(T)\right),\\
\nonumber
& & |\boldsymbol{V}(t_k)|\leq {V}_{max}, \forall t_k\in [0, T],\\
\nonumber
& & \boldsymbol{x}^{(i)}(t_k)\in\mathcal{X}, \forall t_k\in [0, T],
\end{IEEEeqnarray}
We have one remaining challenge, namely the evaluation of the potential force $F(\boldsymbol{x},\boldsymbol{V})$. We recall that $U(\boldsymbol{x}) = a E[(V(\boldsymbol{Y})-\bar{v}(\boldsymbol{X}))^2|\boldsymbol{X}=\boldsymbol{x}]$, where $\bar{v}(\boldsymbol{x})=E[V(\boldsymbol{Y})|\boldsymbol{X}=\boldsymbol{x}]$, and the expectations are computed over the distribution $\mathcal{N}(\boldsymbol{x},\sigma^2I)$, with $\sigma$ a parameter determined from the geometric properties of the spherical particle.

Given the fixed positions of the electrodes $\boldsymbol{y}_{i,j} = (y_1^{i,j}, y_2^{i,j})$, we can approximate the average potential $\bar{v}(\boldsymbol{x})$ and the potential energy, by approximating the conditional Gaussian distribution with a probability mass function. It follows that $$\bar{v}(\boldsymbol{x})\approx \sum_{i,j}p_{i,j}V_{i,j},$$
where $p_{i,j}$ are, not surprisingly, the normalized 2D capacitances in (\ref{eq:03250923}).
Similarly, we have that
$$U(\boldsymbol{x})\approx \sum_{i,j}p_{i,j}\left[V_{i,j}-\bar{v}(\boldsymbol{x})\right]^2$$

The challenge with this approach is that the accuracy of the approximation depends on the granularity of the electrode array: the more electrodes we have, the better the approximations are. We can circumvent this dependency by looking at the electrode potentials as a continuum rather than discrete points. In such a case, the discrete electric potentials become evaluations of a continuous map of potentials $\mathcal{V}(\boldsymbol{y},t)$, evaluated at discrete points, i.e., $V_{i,j}(t) = \mathcal{V}(\boldsymbol{y}_{i,j},t)$. Thus, instead of solving (\ref{eqn:03242113}) to generate a set of discrete potentials over time, we can learn a continuous, parameterized map $\mathcal{V}(\boldsymbol{y},t;\beta)$, where $\beta$ are the parameters of the map.

One additional consequence of learning a continuous map is that we are no longer bound to use the discrete electrode positions when computing the potential energy $U$ and its gradient. We now have the option to use discretization schemes that offer a better accuracy when computing the expectations. In particular, we make use of Gauss quadrature rules \cite{Golub:1967:CGQ:891662}, often found in the theory of generalized chaos polynomials (GPC) \cite{ohagan13-polyn-chaos,Smith:2013:UQT:2568154,citeulike:2621074,Xiu2003}. Gauss quadrature rules  provide the means to accurately evaluate the conditional expectations, using a small number of points. Since the conditional probability distribution of $\boldsymbol{Y}|\boldsymbol{X}=\boldsymbol{x}$ can be expressed as product of two Gaussian distributions $Y_i|X_i=x_i\backsim \mathcal{N}(x_i,\sigma^2)$, with $i\in\{1,2\}$, we have that the expectation of $V(\boldsymbol{Y})|\boldsymbol{X}=\boldsymbol{x}$ is given by
$$\bar{v}(\boldsymbol{x}) = E[V(\boldsymbol{Y})|\boldsymbol{X}=\boldsymbol{x}]$$
\begin{equation}
\label{eqn:03292142}
\approx \frac{1}{{\pi}}\sum_{i,j=1}^nw_iw_jV\left(\sqrt{2}\sigma y_i+x_1,\sqrt{2}\sigma y_j+x_2\right),
\end{equation}
where $n$ is the number of sample points, $y_i$ are the roots of the physicists' version of the Hermite polynomial $H_n(y)$ and $w_i$ are associated weights given by $w_i = \frac{2^{n-1}n!\sqrt{\pi}}{n^2[H_{n-1}(y_i)]^2}$. Similarly, the variance of $V(\boldsymbol{Y})|\boldsymbol{X}=\boldsymbol{x}$ can be approximated as
$$U(\boldsymbol{x}) = E\left[\left(V(\boldsymbol{Y})-\bar{v}(\boldsymbol{X})\right)^2|\boldsymbol{X}=\boldsymbol{x}\right]$$
\begin{equation}
\label{eqn:03292143}
\approx\frac{1}{{\pi}}\sum_{i,j=1}^nw_iw_j\left[V\left(\sqrt{2}\sigma y_i+x_1,\sqrt{2}\sigma y_j+x_2\right)-\bar{v}(\boldsymbol{x})\right]^2.
\end{equation}

We now have the optimal control formulation for learning a continuous map of potentials over time:
\begin{IEEEeqnarray}{rcl}
\label{eqn:03242135}
\min_{\beta,\boldsymbol{x}^{(i)}(t_k)} & &  \sum_{i,j} (\hat{f}(\boldsymbol{x}_{i,j}, T)-f_d(\boldsymbol{x}_{i,j}))^2\\
\nonumber
\textmd{s.t.: } &  & \boldsymbol{x}_{t_{k+1}}^{(i)} = \boldsymbol{x}_{t_{k}}^{(i)} + \frac{\Delta t}{2}\left[F(\boldsymbol{x}^{(i)}(t_{k+1}), \mathcal{V}(\cdot,t_{k+1};\beta))\right.\\
\nonumber
& &\left. +F(\boldsymbol{x}^{(i)}(t_{k}), \mathcal{V}(\cdot,t_{k};\beta))\right], \boldsymbol{x}^{(i)}(0)=\boldsymbol{x}^{(i)}_{0}\\
\nonumber
& & \hat{f}(\boldsymbol{x}_{i,j},T) = \frac{1}{n}\sum_{i=1}^n K_H\left(\boldsymbol{x}_{i,j} - \boldsymbol{x}^{(i)}(T)\right),\\
\nonumber
& & |\mathcal{V}(\cdot,t_{k+1};\beta)|\leq {V}_{max}, \forall t_k\in [0, T],\\
\nonumber
& & \boldsymbol{x}^{(i)}(t_k)\in\mathcal{X}, \forall t_k\in [0, T].
\end{IEEEeqnarray}

The advantage of learning a control map is that we can control the complexity of the optimization variables as the number of electrodes increase. There remains the challenge of selecting a parameterization for $ \mathcal{V}(\boldsymbol{x},t;\beta)$. We can take a global approach to the map representation and get inspiration from spectral methods \cite{Trefethen96finitedifference} to represent the map $\mathcal{V}: \mathcal{X}\times [0, T]$ in terms of a set of polynomial basis functions, namely:
$$\mathcal{V}(\boldsymbol{x},t) = \sum_{i,j=1}^Mv_{i,j}(t)\varphi_{i,j}(\boldsymbol{x}),$$
where $M^2$ is the number of terms in the approximation, and $\varphi_{i,j}(\boldsymbol{x}) = \varphi_{i}(x_1)\varphi_{i}(x_2)$, with $\{\varphi_{i}\}$ a set of polynomial basis functions (e.g., Chebyshev or Legendre polynomials). In this case the parameters $\beta$ are the coefficients $v_{i,j}(t_k)$ for all discrete time instances. Alternatively, we can use universal function approximators, such as neural networks (NNs) \cite{HORNIK1989359}, and the parameters $\beta$ are the weights and biases, i.e., $\mathcal{V}(\boldsymbol{x},t) = \mathcal{V}(\boldsymbol{x},t;\beta)$. This option has the advantage that $\beta$ no longer depends directly on the number of time instances.   Similar ideas can be used to represent the trajectory of the particles. For example, we can use NNs that take as input time and generate as output the particle position. Such approaches have already been used in the context of PDEs, where NNs are used to approximate PDE solutions \cite{https://doi.org/10.1002/gamm.202100006}. Unlike more traditional methods though (e.g., finite elements, (pseudo-)spectral methods), the effects of the approximation errors are much more difficult to quantify.

\begin{remark}
  The loss function in the optimal control formulation includes the KDE for the particle density at the final time only. It is well understood that the quality of the estimator depends on the choice of the bandwidth matrix $H$. At the final time, we can use the statistics of target density to select $H$. For example, using Silverman's rule (Scott's rule is identical for the 2D case),  we can choose $\sqrt{H}_{ii}=n^{-1/5}\sigma_i$, and $H_{ij} = 0$, where $\sigma_i$ is the standard deviations for the $i^{th}$ variable, and $i,j\in \{1,2\}$. Since we select the target density, $\sigma_i$ can picked in relation with this density. For the intermediate particle densities, we do not have a good way to estimate $\sigma_i$, hence we may end up computing under or over smooth estimates.
\end{remark}

\section{Results}
\label{sec:Results}
We use global parameterizations for both the particle trajectories and electrodes electric potentials. In particular, the vector of positions $\boldsymbol{Z}:[0,T]\rightarrow \mathds{R}^{2n}$, with $\boldsymbol{Z}=[\boldsymbol{x}^{(1)}, \ldots, \boldsymbol{x}^{(n)}]$, and the electric potentials $V:\mathds{R}^2\times [0,T]\rightarrow \mathds{R}$ are defined as NNs. One advantage of this type of parameterization is that we can use batch execution to evaluate the particle positions and the electric potential at a sequence of time instances and positions, jointly. Another advantage is that there is no longer the need to explicitly discretize the particle dynamics since we can use AD to evaluate the time derivatives. To ensure scalability with the number of optimization variables, we use Adam \cite{Adam}, a first order gradient-based algorithm, and recast the optimization problem in primal-dual flavor. We minimize the loss function:
$$ \sum_{i,j} \left(\frac{1}{n}\sum_{l=1}^n K_H\left(\boldsymbol{x}_{i,j} - \boldsymbol{x}^{(l)}(T)\right)-f_d(\boldsymbol{x}_{i,j})\right)^2 + $$
$${\lambda}\sum_{l=1}^n\sum_{k,i,j}\left\|\dot{\boldsymbol{x}}_{t_{k}}^{(l)} - F(\boldsymbol{x}^{(l)}(t_{k}), \mathcal{V}(\boldsymbol{x}_{i,j},t_{k}))\right\|^2,$$
in terms of the weights and biases of the NNs, using Adam. The bounds constraints on the particle positions and electric potentials are imposed through projections. Periodically, the weight of the constraint $\lambda$ is updated to reflect how far we are from satisfying the constraint, using a projected gradient step:
$$\lambda \leftarrow \left[\lambda + \alpha \left(\sum_{l=1}^n\sum_{k,i,j}\left\|\dot{\boldsymbol{x}}_{t_{k}}^{(l)} - F(\boldsymbol{x}^{(l)}(t_{k}), \mathcal{V}(\boldsymbol{x}_{i,j},t_{k}))\right\|^2 - \epsilon\right)\right]_{+},$$
for a positive stepsize $\alpha$ and a small positive scalar $\varepsilon$, playing the role of tolerance.

\begin{figure}[!htp]
  \begin{center}
    \includegraphics[width=0.7\textwidth]{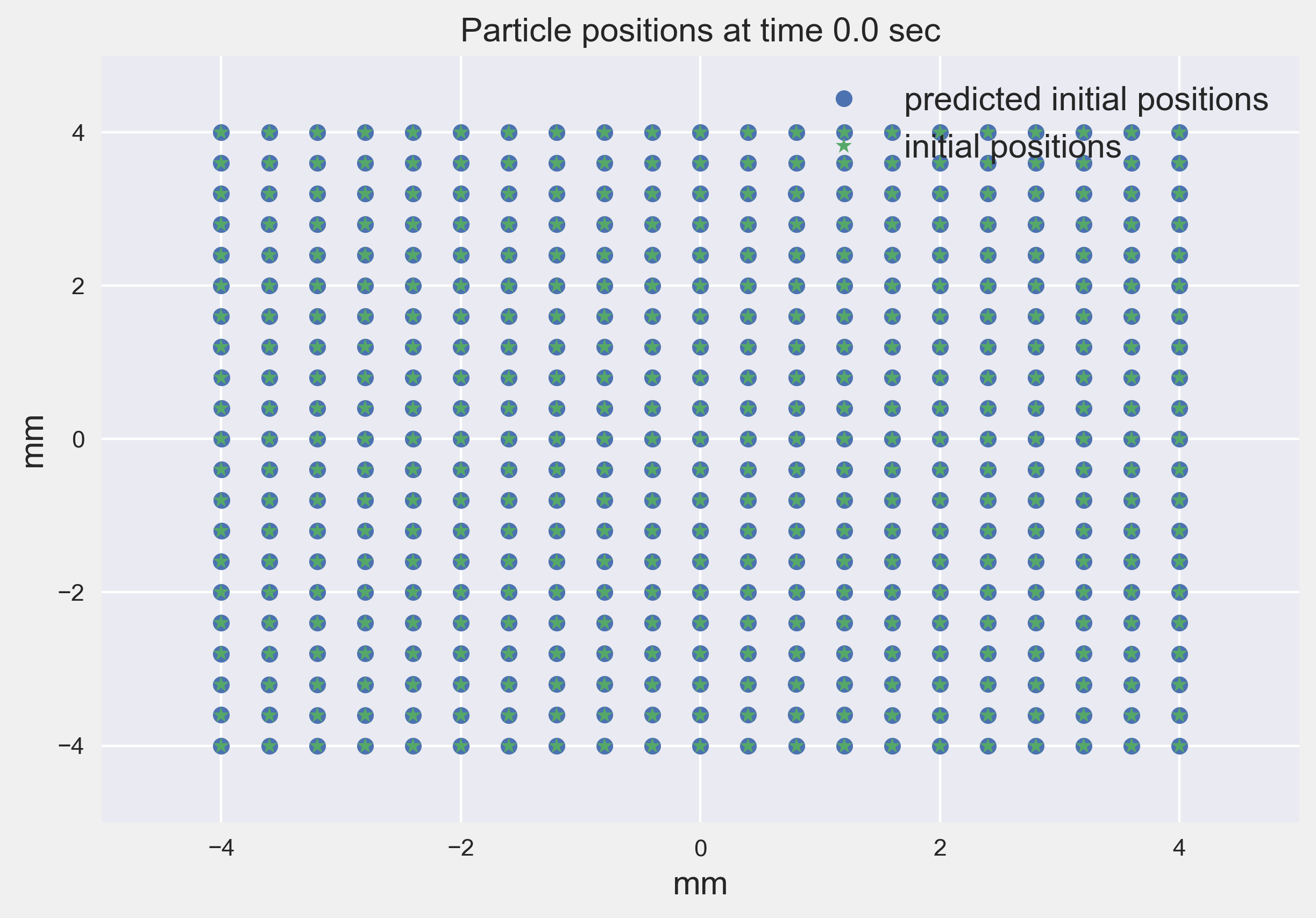}
    \end{center}
  \caption{Particle initial positions.}
  \label{fig:initial}
\end{figure}

\begin{figure}[!htp]
  \begin{center}
    \includegraphics[width=0.7\textwidth]{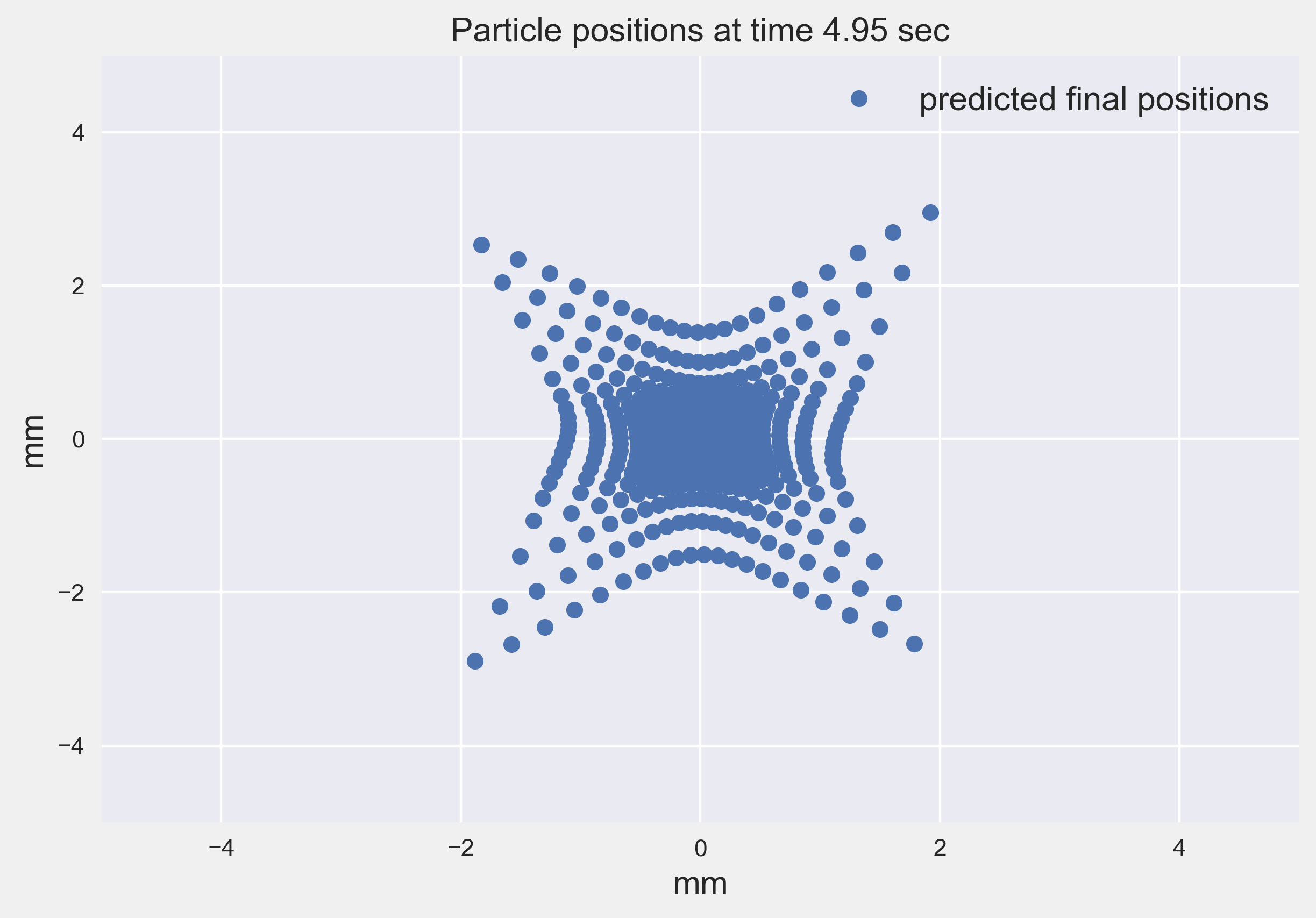}
    \end{center}
  \caption{Particle final positions.}
  \label{fig:final}
\end{figure}

\begin{figure}[!htp]
  \begin{center}
    \includegraphics[width=0.7\textwidth]{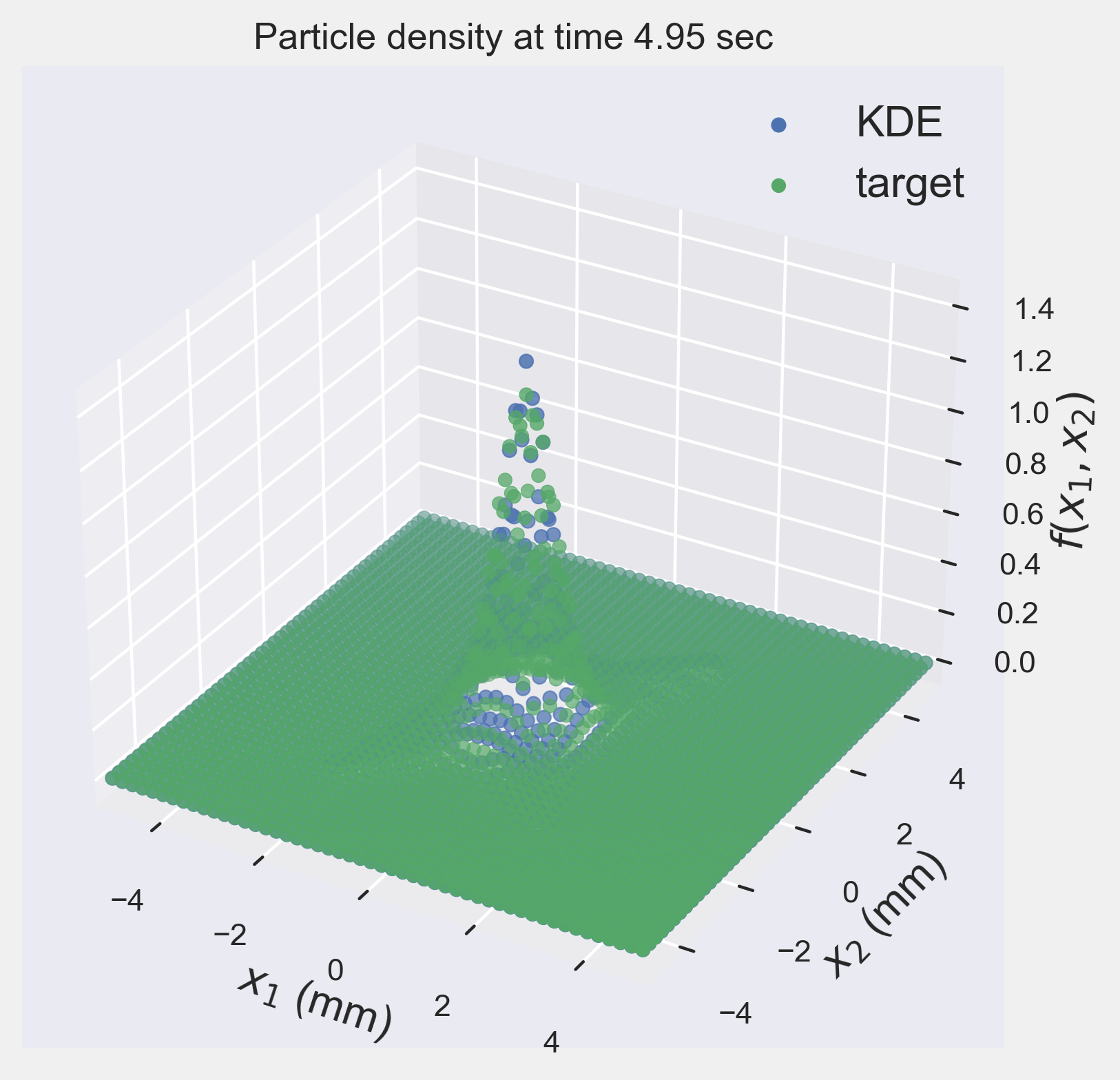}
    \end{center}
  \caption{Particle density comparison: KDE vs target.}
  \label{fig:particle density comparison}
\end{figure}

In our simulation results, we consider $\mathcal{X} = 5 mm\times 5mm$, and the electrodes are uniformly distributed, with $0.25$ mm electrode pitch, resulting in 1680 electrodes. The maximum electrode potential magnitude is 400V. The particle capacitance is a multivariate, Gaussian distribution with mean at the particle position $\boldsymbol{x}$ and covariance matrix $\sigma^2 I$, where $\sigma=0.25\times 10^{-3}$. We consider a time horizon of 5 seconds and a sampling period of $0.05$ msec. We assume that the camera is able to identify the positions of 450 particles and these positions will be used as initial conditions in the optimal control problem. We use two NNs to model  each of the $x_1$ and $x_2$ directions of the particle positions. These NNs have one hidden layer of size 1500, using {\tt tanh} as activation function. The electric potential map $\mathcal{V}$ is modeled as a NN with a hidden layer of size 500 and {\tt tanh} as activation function. Note that while the complexity of the NN modeling the map $\mathcal{V}$ can remain constant, the complexity of the NNs modeling the particle positions over time will depend on the number of particles. We set the target particle distribution to a multivariate Gaussian distribution centered at zero, with covariance matrix $\sigma^2 I$, where $\sigma=0.5\times 10^{-3}$. The bandwidth parameter is chosen as $h=1.06\sigma n^{-1/5}$, where $n$ is the number of particles. In addition, we set the tolerance for satisfying the particle dynamics constraint to $\varepsilon=10^{-4}$. Note that this tolerance is significant considering that it is applied to the sum (not the average) of more than $10^7$ residuals, when accounting for the number of particles, time samples, and grid points at which the electric potentials are evaluated. We start from a perfect uniform particle distribution, as shown in Figure \ref{fig:initial}. After solving the optimal control problem, the particle positions after 5 seconds are shown in Figure \ref{fig:final}. The positions are computed by evaluating the trained NNs at the end of the time horizon. In Figure \ref{fig:initial} we plot both the initial conditions predicted by the NNs and the true initial conditions, showing that the NNs perfectly match the ground truth. Figure \ref{fig:particle density comparison} depicts a comparison between the target density and the KDE computed using the particle positions, at the end of the time horizon. The MSE between the KDE and the target density is $6.7\times 10^{-4}$. It is well in the range of errors obtained by sampling from the Gaussian distribution and computing the KDE using these samples.

\begin{figure}
\begin{subfigure}{.5\textwidth}
  \centering
  \includegraphics[width=1\textwidth]{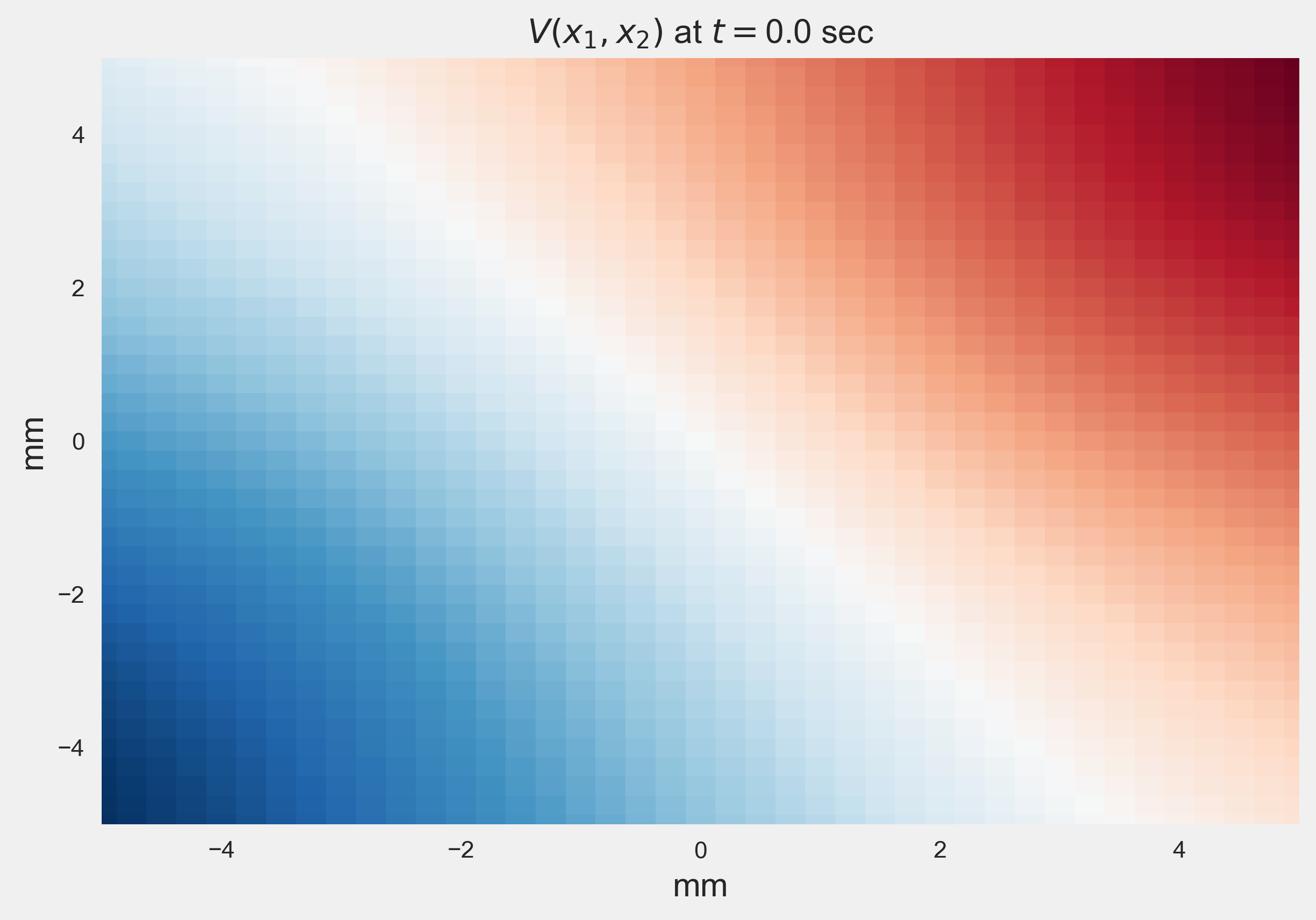}
  \caption{t=0 sec.}
  \label{fig:potentials_t0}
\end{subfigure}%
\begin{subfigure}{.5\textwidth}
  \centering
  \includegraphics[width=1\textwidth]{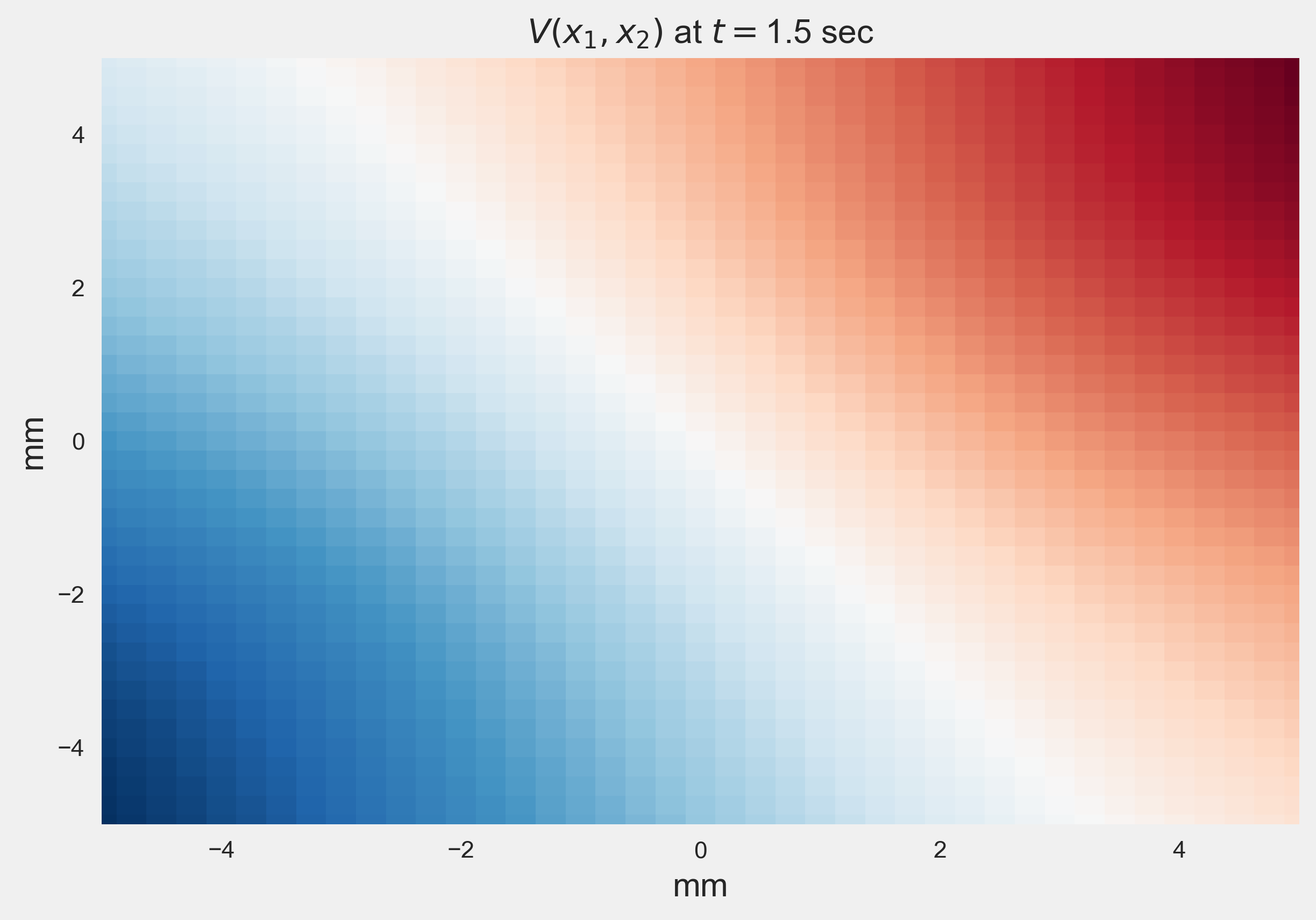}
  \caption{t=1.5 sec.}
  \label{fig:potentials_t1}
\end{subfigure}\\
\begin{subfigure}{.5\textwidth}
  \centering
  \includegraphics[width=1\textwidth]{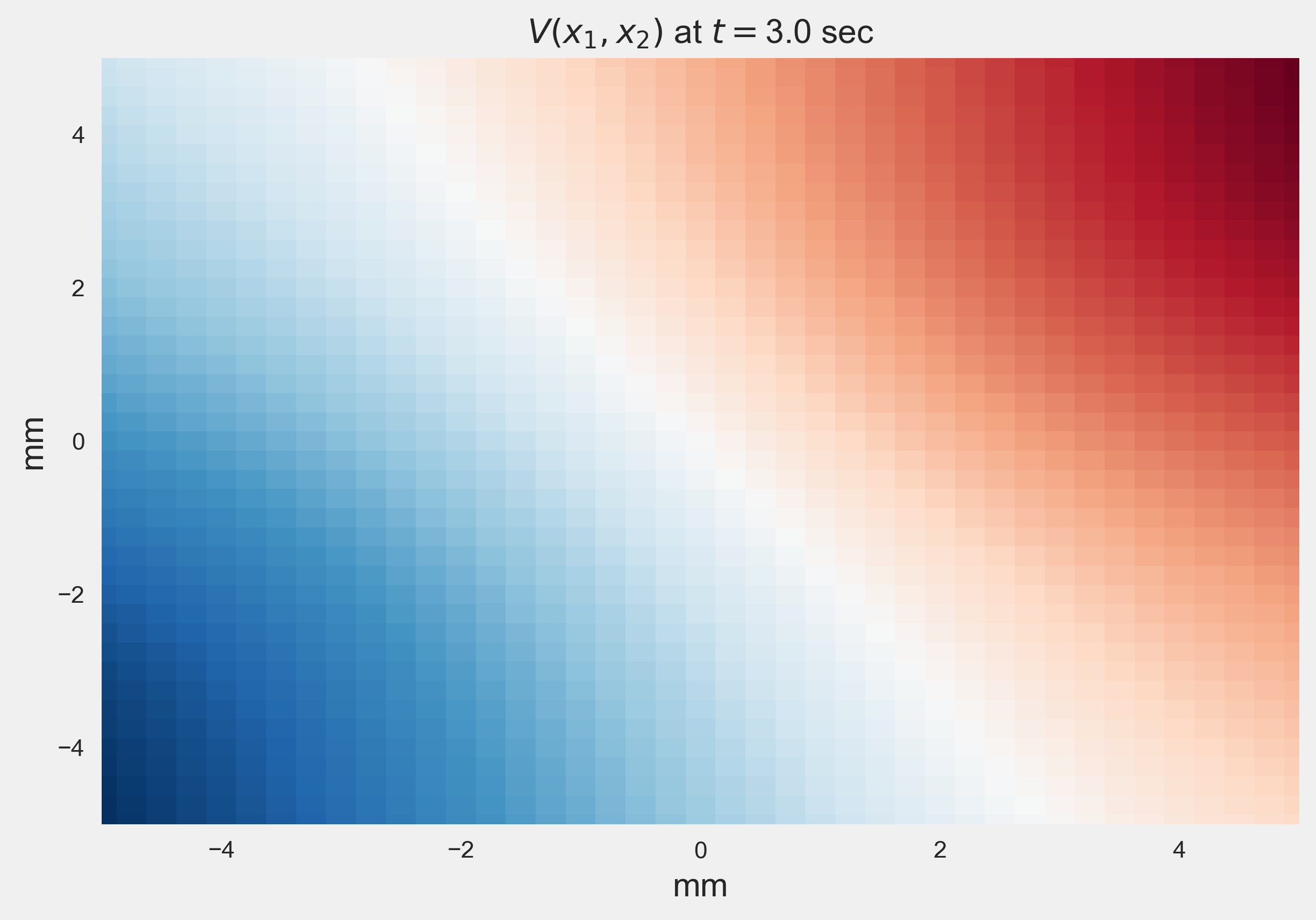}
  \caption{t=3 sec.}
  \label{fig:potentials_t2}
\end{subfigure}%
\begin{subfigure}{.5\textwidth}
  \centering
  \includegraphics[width=1\textwidth]{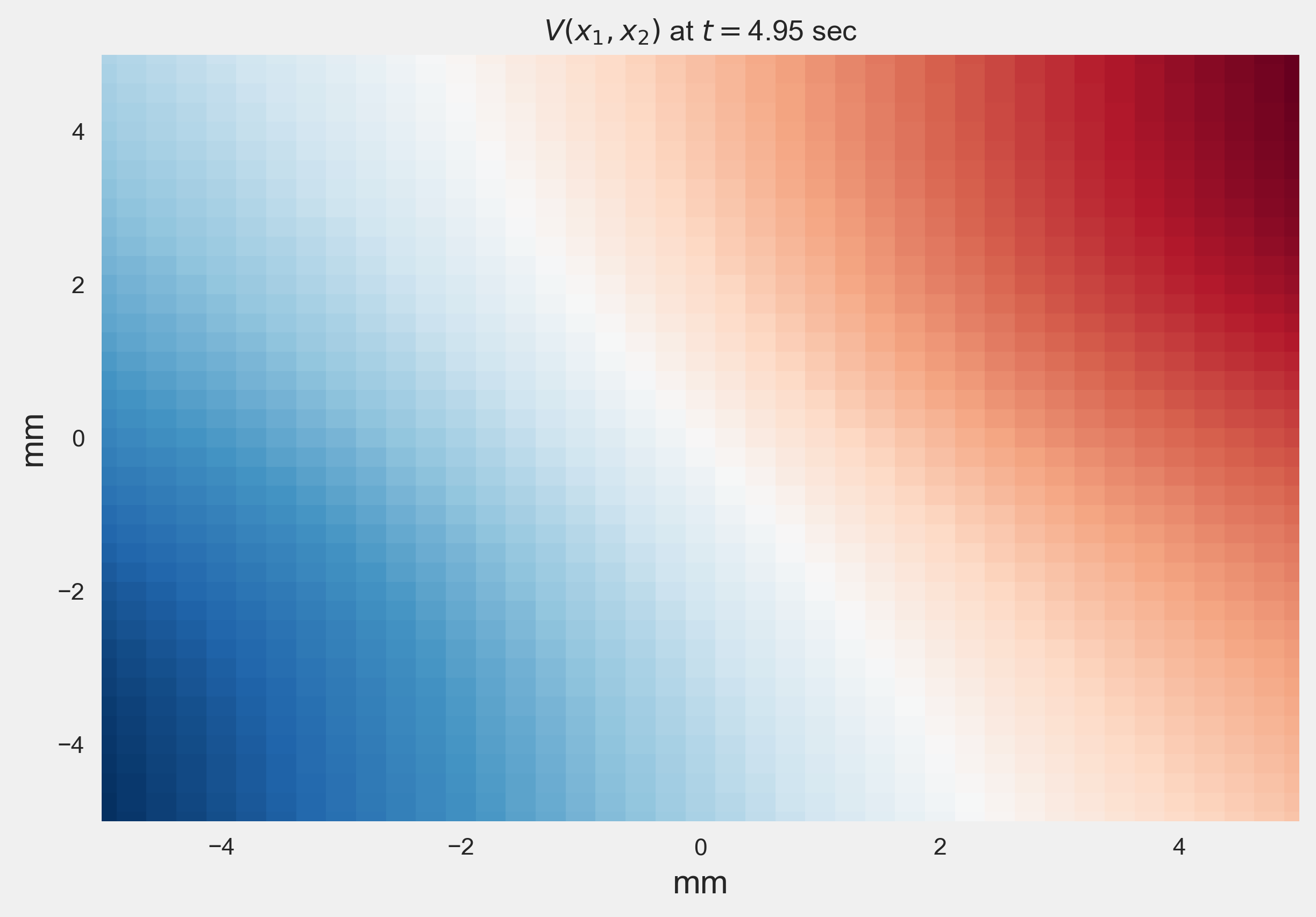}
  \caption{t=4.95 sec.}
  \label{fig:potentials_t4}
\end{subfigure}
\caption{Electrode electric potentials over time.}
\label{fig:electric potentials}
\end{figure}

Figure \ref{fig:electric potentials} shows the electric potential changes over the time horizon, for four time samples. The red and blue colors represent negative and positive potentials, respectively. The intensity of the colors is proportional to the magnitude of the potentials. Interestingly, the plots show a clockwise rotation of the potential distribution over time. The results for the electric potentials are not unique and depend on the type of parameterization we use. As we increase the complexity of the NN, we can model a larger class of control inputs while increasing the complexity of the optimization problem. The NN modeling the behavior of the electrode potentials over time has a total of 2501 parameters. Had we considered the control inputs at each electrode and each time sample, we would have had more than $1.6\times 10^4$ optimization variables just for the control input. The complexity reduction by using a NN (or other type of parameterization) is obvious. What is not obvious upfront is what NN architecture to start with. In our case, we used a parsimonious approach: we started with a simple, one hidden layer architecture and increased its complexity (i.e., the size of the hidden layer) until we obtain a satisfactory result. The same approach was used for modeling the particle positions.

\section{Conclusions}
We addressed the problem of shaping the distribution of particles immersed in a dielectric fluid, by manipulating an electric field controlled by an array of electrodes. We employed a KDE to approximate the particle density, where the dynamics of a particle was determined  using a 2D capacitive-based model of motion. We provided a probabilistic view for interpreting  the particle dynamics for an arbitrarily large number of electrodes. In addition, we showed how we can use Gauss-Hermite quadrature to accurately approximate the potential energy of the particle. We formulated an optimal control problem that minimizes the $L_2$ norm between the particle density at the end of a time horizon and a target density, having the particle dynamics as constraint. We used automatic differentiation to compute derivatives of physical quantities (e.g., potential energy) and the gradient of the cost and constraint functions. We demonstrated our approach by shaping the density of particles from a uniform to a Gaussian distribution. As future work, we will compare the KDE-based approach to an optimal control formulation that uses the Liouville equation, as dynamical constraint.
\bibliographystyle{plain}
\bibliography{references,refs}

\end{document}